\documentclass{elsart}

\usepackage{graphicx}
\usepackage{amssymb,amsmath,amsfonts,latexsym,euscript}

\usepackage{amssymb}

\def\d{\,{\rm d}}

\def\abd{\mathop {\rm \leftharpoondown \!\!\!\!\! \rightharpoondown} }

\renewcommand{\theequation}{\thesection.\arabic{equation}}
\newtheorem{theorem}[subsection]{Theorem}
\newtheorem{lemma}[subsection]{Lemma}

\newtheorem{conjecture}[subsection]{Conjecture}

\begin{document}

\begin{frontmatter}



\title{The discrete analogue of the differential operator
$\frac{\d^{2m}}{\d x^{2m}}+2\omega^2\frac{\d^{2m-2}}{\d
x^{2m-2}}+\omega^4\frac{\d^{2m-4}}{\d x^{2m-4}}$}


\author{A.R. Hayotov},
\ead{abdullo\_hayotov@mail.ru}

\address{Institute of Mathematics, National University of Uzbekistan,
Do`rmon yo`li str.,29, Tashkent 100125, Uzbekistan}

\begin{abstract}
In the present paper we construct the discrete analogue
$D_m(h\beta)$ of the differential operator $\frac{\d^{2m}}{\d
x^{2m}}+2\omega^2\frac{\d^{2m-2}}{\d
x^{2m-2}}+\omega^4\frac{\d^{2m-4}}{\d x^{2m-4}}$. The discrete
analogue $D_m(h\beta)$ plays the main role in construction of
optimal quadrature formulas and interpolation splines minimizing
the semi-norm in the $K_2(P_m)$ Hilbert space.
\end{abstract}

\end{frontmatter}

\section{Introduction and Preliminaries}

The optimization problem of approximate integration formulas in
the modern sense appears as the problem of finding the minimum of
the norm of a error functional $\ell$ given on some set of
functions.

The minimization problem of the norm of the error functional by
coefficients was reduced in \cite{Sob74,SobVas} to the system of
difference equations of Wiener-Hopf type in the space $L_2^{(m)}$,
where $L_2^{(m)}$ is the space of functions with square integrable
$m-$th generalized derivative. Existence and unique\-ness of the
solution of this system was proved by S.L. Sobolev. In the works
\cite{Sob74,SobVas} the description of some analytic algorithm for
finding the coefficients of optimal formulas is given. For this
S.L. Sobolev defined and investigated the discrete analogue
$D_{hH}^{(m)}(h\beta)$ of the poly\-har\-monic  operator
$\Delta^m$. The problem of construction of the discrete opera\-tor
$D_{hH}^{(m)}(h\beta)$ for $n-$ dimensional case was very hard. In
one dimensional case the discrete analogue $D_{h}^{(m)}(h\beta)$
of the differential operator $\frac{\d^{2m}}{\d x^{2m}}$ was
constructed by Z.Zh. Zhamalov \cite{Zham78} and Kh.M. Shadimetov
\cite{Shad85}.

Further, in the work \cite{ShadHay04_1} the discrete analogue of
the differential operator $\frac{\d^{2m}}{\d
x^{2m}}-\frac{\d^{2m-2}}{\d x^{2m-2}}$ was constructed. The
constructed discrete analogue of the operator $\frac{\d^{2m}}{\d
x^{2m}}-\frac{\d^{2m-2}}{\d x^{2m-2}}$ was applied for finding the
coefficients of the optimal quadrature formulas (see
\cite{ShadHay_04_2,ShadHay13}) and for construction of
interpolation splines minimizing the semi-norm (see
\cite{ShadHay12}) in the space $W_2^{(m,m-1)}(0,1)$, where
$W_2^{(m,m-1)}(0,1)$ is the Hilbert space of functions $\varphi$
which $\varphi^{(m-1)}$ is absolutely continuous, $\varphi^{(m)}$
belongs to $L_2(0,1)$ and
$\int_0^1(\varphi^{(m)}(x)+\varphi^{(m-1)}(x))^2\d x<\infty$.

Here, we  mainly use a concept of functions of a discrete argument
and the corresponding operations \cite[Chapter VII]{Sob74}.  For
completeness we give some of definitions.

Let  $\beta\in \mathbb{Z}$, $h=\frac{1}{N}$, $N=1,2,...$ . Assume
that $\varphi$ and $\psi$ are real-valued functions defined on the
real line $\mathbb{R}$.
\smallskip

{\sc Definition 1.} The function $\varphi(h\beta)$ is a
\emph{function of discrete argument} if it is defined on some set
of integer values of $\beta $.
\smallskip

 {\sc Definition 2.} The \emph{inner product} of two
discrete argument functions $\varphi(h\beta)$ and $\psi(h\beta)$
is defined as
\[
\left[ {\varphi,\psi} \right] = \sum_{\beta  =  - \infty }^\infty
\varphi(h\beta) \cdot \psi(h\beta).\] Here we assume that the
series on the right hand side  converges absolutely.
\smallskip

{\sc Definition 3.} The \textit{convolution} of two functions
$\varphi(h\beta)$ and $\psi(h\beta)$ is defined as the inner
product
\[
\varphi(h\beta)*\psi(h\beta) = \left[
{\varphi(h\gamma),\psi(h\beta-h\gamma)} \right] = \sum_{\gamma  =
- \infty }^\infty  {\varphi (h\gamma)\cdot \psi(h\beta  -
h\gamma)}.
\]

{\sc Definition 4.} The function
$\stackrel{\abd}{\varphi}\!\!(x)=\sum\limits_{\beta=-\infty}^{\infty}\varphi(h\beta)\delta(x-h\beta)$
is called \emph{the harrow-shaped function} corresponding to the
function of discrete argument $\varphi(h\beta)$, where $\delta$ is
Dirac's delta function.

We are interested in to find a discrete function $D_m(h\beta)$
that satisfies the following equation
\begin{equation}\label{eq.(1.1)}
D_m(h\beta)*G_m(h\beta)=\delta_{\d}(h\beta),
\end{equation}
where
\begin{eqnarray}
G_m(h\beta)&=&\frac{(-1)^m\mathrm{sign}(h\beta)}{4\omega^{2m-1}}\Bigg[(2m-3)\sin(h\omega\beta)-
h\omega\beta\cos(h\omega\beta)\nonumber \\
&&\qquad \qquad \qquad
+2\sum\limits_{k=1}^{m-2}\frac{(-1)^{k}(m-k-1)(h\omega\beta)^{2k-1}}{(2k-1)!}\Bigg],\label{eq.(1.2)}
\end{eqnarray}
$m\geq 2$, $\omega>0$, $\delta_{\d}(h\beta)$ is equal to 1 when
$\beta=0$ and is equal to 0 when $\beta\neq 0$, i.e.
$\delta_{\d}(h\beta)$ is the discrete delta function.

The discrete function $D_m(h\beta)$ plays an important role in the
calculation of the coefficients of the optimal quadrature formulas
and interpolation splines minimizing a semi-norm in the Hilbert
space
\[
K_2(P_m)=\Bigl\{\varphi:[0,1]\to \mathbb{R}\ \Bigm|\
\varphi^{(m-1)} \mbox{ is abs. cont. and } \varphi^{(m)}\in
L_2(0,1)\Bigr\},
\]
equipped with the norm
\begin{equation}\label{eq.(1.3)}
 \|\varphi\|=\left\{\int\limits_0^1\left(P_m\left(\frac{\d}{\d
x}\right) \varphi(x)\right)^2\d x\right\}^{\frac12},
\end{equation}
where $P_m(\frac{\d}{\d x})=\frac{\d^m}{\d
x^m}+\omega^2\frac{\d^{m-2}}{\d x^{m-2}}$, $\omega>0$ and $
{\int_0^1 {\left( {P_m \left( {{\d \over {\d x}}} \right)\varphi
(x)} \right)} ^2 \d x}<\infty.$

Note that the equality (\ref{eq.(1.3)}) is semi-norm, moreover
$\|\varphi\|=0$ if and only if $\varphi(x)=c_1\sin \omega
x+c_2\cos \omega x+Q_{m-3}(x)$, with $Q_{m-3}(x)$ a polynomial of
degree $m-3$. In particular, the optimal quadrature formulas and
interpolation splines in the $K_2(P_m)$ space are exact for the
trigonometric functions $\sin\omega x$, $\cos \omega x$ and
polynomials of degree $m-3$.

It should be noted that for a linear differential operator of
order $m$,\linebreak $$ L:= P_m(\d/\d x)=\frac{\d^m}{\d
x^m}+a_{m-1}(x)\frac{\d^{m-1}}{\d x^{m-1}}+...+a_1(x)\frac{\d}{\d
x}+a_0(x),$$ Ahlberg, Nilson, and Walsh in the book \cite[Chapter
6]{Ahlb67} investigated the Hilbert spaces in the context of
generalized splines. Namely, with the inner product
\[\langle\varphi,\psi\rangle=\int_0^1 L\varphi (x)\cdot L\psi (x)  \d x ,\]
$K_2(P_m)$ is a Hilbert space if we identify functions that differ
by a solution of $L\varphi=0$.

Note that the equation (\ref{eq.(1.1)}) is the discrete analogue
to the following equation
\begin{equation}\label{eq.(1.4)}
\left(\frac{\d^{2m}}{\d x^{2m}}+2\omega^2\frac{\d^{2m-2}}{\d
x^{2m-2}}+\omega^4\frac{\d^{2m-4}}{\d
x^{2m-4}}\right)G_m(x)=\delta(x),
\end{equation}
where
\begin{eqnarray*}
G_m(x)&=&\frac{(-1)^m\mathrm{sign}(x)}{4\omega^{2m-1}}\Bigg[(2m-3)\sin
(\omega x)-\omega x\cdot \cos(\omega x)\\
&& \qquad\qquad\qquad
+2\sum_{k=1}^{m-2}\frac{(-1)^{k}(m-k-1)(\omega
x)^{2k-1}}{(2k-1)!}\Bigg],
\end{eqnarray*}
and $\delta$ is the Dirac's delta function.

Moreover the discrete function $D_m(h\beta)$ has analogous
properties the differential operator $\frac{\d^{2m}}{\d
x^{2m}}+2\omega^2\frac{\d^{2m-2}}{\d
x^{2m-2}}+\omega^4\frac{\d^{2m-4}}{\d x^{2m-4}}$: the zeros of the
discrete operator $D_m(h\beta)$ coincides with the discrete
functions corresponding to the zeros of the operator
$\frac{\d^{2m}}{\d x^{2m}}+2\omega^2\frac{\d^{2m-2}}{\d
x^{2m-2}}+\omega^4\frac{\d^{2m-4}}{\d x^{2m-4}}$. The discrete
function $D_m(h\beta)$ is called \emph{the discrete analogue}  of
the differential operator $\frac{\d^{2m}}{\d
x^{2m}}+2\omega^2\frac{\d^{2m-2}}{\d
x^{2m-2}}+\omega^4\frac{\d^{2m-4}}{\d x^{2m-4}}$.

The rest of the paper is organized as follows: in Section 2 we
give some known formulas and  auxiliary results which will be used
in the construction of the discrete function $D_m(h\beta)$.
Section 3 is devoted to the construction of  the discrete analogue
$D_m(h\beta)$ of the differential operator $\frac{\d^{2m}}{\d
x^{2m}}+2\omega^2\frac{\d^{2m-2}}{\d
x^{2m-2}}+\omega^4\frac{\d^{2m-4}}{\d x^{2m-4}}$.

\section{Auxiliary results}
\setcounter{equation}{0}

In this section we explain some known formulas (see, for instance,
\cite{Sob74,Vlad79}) and auxiliary results which we use in the
construction of the discrete analogue $D_m(h\beta)$ of the
differential operator $\frac{\d^{2m}}{\d
x^{2m}}+2\omega^2\frac{\d^{2m-2}}{\d
x^{2m-2}}+\omega^4\frac{\d^{2m-4}}{\d x^{2m-4}}$.

For the Fourier transformations and their properties:
\begin{eqnarray}
&&F[\varphi(x)]=\int\limits_{-\infty}^{\infty}\varphi(x)e^{2\pi
ipx}\d x,\ \ \
F^{-1}[\varphi(p)]=\int\limits_{-\infty}^{\infty}\varphi(p)e^{-2\pi
ipx}\d p,\label{eq.(2.1)}\\
&&F[\varphi*\psi]=F[\varphi]\cdot F[\psi], \label{eq.(2.2)}\\
&& F[\varphi\cdot \psi]=F[\varphi]* F[\psi], \label{eq.(2.3)}\\
&& F[\delta^{(\alpha)}(x)]=(-2\pi ip)^{\alpha},\ \ F[\delta(x)]=1,
\label{eq.(2.4)}
\end{eqnarray}
where $*$ is the convolution and for two functions $\varphi$ and
$\psi$ is  defined as follows
$$
(\varphi*\psi)(x)=\int\limits_{-\infty}^{\infty}\varphi(x-y)\psi(y)\d
y=\int\limits_{-\infty}^{\infty}\varphi(y)\psi(x-y)\d y.
$$
For the Dirac's delta function:
\begin{eqnarray}
&&\delta(hx)=h^{-1}\delta(x),\label{eq.(2.5)}\\
&&\delta(x-a)\cdot f(x)=\delta(x-a)\cdot f(a), \label{eq.(2.6)}\\
&& \delta^{(\alpha)}(x)*f(x)=f^{(\alpha)}(x), \label{eq.(2.7)}
\\
&&\phi_0(x)=\sum\limits_{\beta=-\infty}^{\infty}\delta(x-\beta)=\sum\limits_{\beta=-\infty}^{\infty}e^{2\pi
i x\beta}. \label{eq.(2.8)}
\end{eqnarray}

It is known \cite{Froben,Sob06_1,SobVas} that {\it the
Euler-Frobenius polynomials} $E_k (x)$  are defined by the
following formula
\begin{equation}\label{eq.(2.9)}
 E_k (x) = \frac{{(1 - x)^{k + 2} }}{x}\left( {x\frac{d}{{dx}}}
\right)^k \frac{x}{{(1 - x)^2 }},\ \ k=1,2,...
\end{equation}
$E_0 (x) = 1$, moreover all roots $x_j^{(k)}$ of the
Euler-Frobenius polynomial $E_k(x)$ are real, negative and simple,
i.e.
\begin{equation}\label{eq.(2.10)}
x_1^{(k)}<x_2^{(k)}<...<x_k^{(k)}<0.
\end{equation}
Furthermore, the roots equally spaced from the ends of the chain
(\ref{eq.(2.10)}) are reciprocal, i.e.
$$
x_j^{(k)}\cdot x_{k+1-j}^{(k)}=1.
$$
Euler obtained the following formula for the coefficients
$a_s^{(k)}$, $s=0,1,...,k$ of the polynomial
$E_k(x)=\sum\limits_{s=0}^ka_s^{(k)}x^s$:
\begin{equation}\label{eq.(2.11)}
a_s^{(k)}=\sum\limits_{j=0}^s(-1)^j{k+2\choose j} (s+1-j)^{k+1}.
\end{equation}

The polynomial $E_k(x)$ satisfies the following identity
\begin{equation}\label{eq.(2.12)}
E_k(x)=x^kE_k\left(\frac{1}{x}\right),\ \ x\neq 0,
\end{equation}
i.e. $a_s^{(k)}=a_{k-s}^{(k)},\ \ s=0,1,2,...,k$.

Now we consider the following polynomial of degree $2m-2$
$$
{\cal P}_{2m-2}(x)=\sum_{s=0}^{2m-2}p_s^{(2m-2)}x^s
=(1-x)^{2m-4}\bigg[[(2m-3)\sin h\omega-h\omega \cos h\omega]x^2
$$
$$
\qquad\qquad\qquad +[2h\omega-(2m-3)\sin(2h\omega)]x+[(2m-3)\sin
h\omega-h\omega \cos h\omega]\bigg]
$$
\begin{equation}\label{eq.(2.13)}
+2(x^2-2x\cos
h\omega+1)^2\sum_{k=1}^{m-2}\frac{(-1)^k(m-k-1)(h\omega)^{2k-1}(1-x)^{2m-2k-4}E_{2k-2}(x)}{(2k-1)!},
\end{equation}
where $E_{2k-2}(x)$ is the Euler-Frobenius polynomial of degree
$2k-2$, $\omega>0$, $h\omega\leq 1$, $h=1/N$, $N\geq m-1$, $m\geq
2$.

The polynomial ${\cal P}_{2m-2}(x)$ has the following properties.

\begin{lemma}\label{L2.1}
The Polynomial ${\cal P}_{2m-2}(x)$ satisfies the identity
$$
{\cal P}_{2m-2}(x)=x^{2m-2}{\cal
P}_{2m-2}\left(\frac{1}{x}\right).
$$
\end{lemma}

\begin{lemma}\label{L2.3}
For the derivative of the polynomial ${\cal P}_{2m-2}(x)$ the
following equality holds
$$
{\cal
P}_{2m-2}'\left(\frac{1}{x_k}\right)=-\frac{1}{x_k^{2m-4}}{\cal
P}_{2m-2}'(x_k),
$$
where $x_k$ is a root of the polynomial ${\cal P}_{2m-2}(x)$.
\end{lemma}

Further we prove Lemmas \ref{L2.1}-\ref{L2.3}.

\emph{Proof of Lemma \ref{L2.1}.} If we replace $x$ by $\frac1x$
in (\ref{eq.(2.13)}), we obtain
$$
{\cal
P}_{2m-2}\left(\frac1x\right)=\left(1-\frac1x\right)^{2m-4}\bigg[[(2m-3)\sin
h\omega-h\omega \cos h\omega]\frac{1}{x^2}
$$
$$
\qquad\qquad\qquad
+[2h\omega-(2m-3)\sin(2h\omega)]\frac1x+[(2m-3)\sin
h\omega-h\omega \cos h\omega]\bigg]
$$
$$
+2\left(\frac{1}{x^2}-\frac{2}{x}\cos
h\omega+1\right)^2\sum_{k=1}^{m-2}\frac{(-1)^k(m-k-1)(h\omega)^{2k-1}\left(1-\frac1x\right)^{2m-2k-4}E_{2k-2}\left(\frac{1}{x}\right)}{(2k-1)!}.
$$
Hence, taking into account (\ref{eq.(2.12)}) and the fact that
$\left(1-\frac{1}{x}\right)^{2k}=\frac{1}{x^{2k}}(1-x)^{2k}$, we
get the statement of the lemma. Lemma \ref{L2.1} is proved.\hfill
$\Box$

\emph{Proof of Lemma \ref{L2.3}.} Let $x_1,x_2,...,x_{2m-2}$ be
the roots of the polynomial ${\cal P}_{2m-2}(x)$. Then from Lemma
\ref{L2.1} we immediately get that if $x_k$ is the root then
$\frac{1}{x_k}$ is also the root of the polynomial ${\cal
P}_{2m-2}(x)$ and if $|x_k|<1$ then $|\frac{1}{x_k}|>1$.

Suppose the roots $x_1,x_2,...,x_{2m-2}$ are situated  such that
\begin{equation}\label{eq.(2.14)}
x_k\cdot x_{2m-1-k}=1.
\end{equation}
Therefore we have
$$
{\cal P}_{2m-2}(x)=p_{2m-2}^{(2m-2)}(x-x_1)(x-x_2)...(x-x_{2m-2}),
$$
where $p_{2m-2}^{(2m-2)}$ is the leading coefficient of the
polynomial ${\cal P}_{2m-2}(x)$.\\
Hence
$$
{\cal
P}_{2m-2}'(x)=p_{2m-2}^{(2m-2)}\sum\limits_{j=1}^{2m-2}\frac{\prod\limits_{i=1}^{2m-2}(x-x_i)}{x-x_j}.
$$
Then for $x=x_k$ we get
\begin{equation}\label{eq.(2.15)}
{\cal P}_{2m-2}'(x_k)=p_{2m-2}^{(2m-2)}\prod\limits_{i=1,i\neq
k}^{2m-2}(x_k-x_i)
\end{equation}
and for $x=\frac{1}{x_k}$, taking into account (\ref{eq.(2.14)}),
using the equality
$$
\prod\limits_{i=1,i\neq
2m-1-k}^{2m-2}\frac{1}{x_{2m-1-i}}=\frac{1}{x_1x_2...x_{k-1}x_{k+1}x_{k+2}...x_{2m-2}}=\frac{1}{x_{2m-1-k}}=x_k
$$
we have
\begin{eqnarray}
 {\cal
P}_{2m-2}'\left(\frac{1}{x_k}\right)&=&p_{2m-2}^{(2m-2)}\sum\limits_{j=1}^{2m-2}\frac{\prod\limits_{i=1}^{2m-2}(\frac{1}{x_k}-x_i)}{\frac{1}{x_k}-x_j}\nonumber\\
&=&p_{2m-2}^{(2m-2)}\frac{1}{x_k^{2m-3}}\prod\limits_{i=1,i\neq
2m-1-k}^{2m-2}(1-x_kx_i)\nonumber\\
&=&-p_{2m-2}^{(2m-2)}\frac{1}{x_k^{2m-3}}\prod\limits_{i=1,i\neq
2m-1-k}^{2m-2}\frac{x_k-x_{2m-1-i}}{x_{2m-1-i}}\nonumber\\
&=&-\frac{1}{x_k^{2m-4}}\ p_{2m-2}^{(2m-2)}\prod\limits_{i=1,i\neq
k}^{2m-2}(x_k-x_i).\label{eq.(2.16)}
\end{eqnarray}
From (\ref{eq.(2.15)}) and  (\ref{eq.(2.16)}) we get the statement
of the lemma. Lemma \ref{L2.3} is proved.\hfill $\Box$

In the proof of the main result we have to calculate the following
series
\begin{eqnarray}
S_1&=&\sum\limits_{\beta=-\infty}^{\infty}\frac{1}{[\beta-h(p+\frac{\omega}{2\pi})][\beta-h(p-\frac{\omega}{2\pi})]},\label{eq.(2.17)}\\
S_2&=&\sum\limits_{\beta=-\infty}^{\infty}\frac{1}{[\beta-h(p+\frac{\omega}{2\pi})]^2},\label{eq.(2.18)}\\
S_3&=&\sum\limits_{\beta=-\infty}^{\infty}\frac{1}{[\beta-h(p-\frac{\omega}{2\pi})]^2},\label{eq.(2.19)}\\
F[h\stackrel{\abd}{G}_{k,1}]&=&\frac{(-1)^k}{(2\pi)^{2k}}\sum\limits_{\beta=-\infty}^{\infty}\frac{1}{(p-h^{-1}\beta)^{2k}}.\label{eq.(2.20)}
\end{eqnarray}

We denote $\lambda=e^{2\pi iph}$ then the following holds

\begin{lemma}\label{L2.4} For the series \emph{(\ref{eq.(2.17)})-(\ref{eq.(2.20)})}
the following are taken place
$$
\begin{array}{rclrcl}
S_1&=&\frac{-(2\pi)^2\lambda\sin h\omega}{h\omega
(\lambda^2-2\lambda\cos h\omega)}, \ \ \  & S_2&=&\frac{-(2\pi)^2\lambda}{(\lambda^2+1)\cos h\omega -2\lambda+i(\lambda^2-1)\sin h\omega},\\
F[h\stackrel{\abd}{G}_{k,1}]&=&\frac{h^{2k}\lambda\
E_{2k-2}(\lambda)}{(2k-1)!(1-\lambda)^{2k}},&
S_3&=&\frac{-(2\pi)^2\lambda}{(\lambda^2+1)\cos h\omega
-2\lambda+i(\lambda^2-1)\sin h\omega},
\end{array}
$$
where $E_{2k-2}(\lambda)$ is the Euler-Frobenius polynomial of
degree $2m-2$.
\end{lemma}

\emph{Proof.} To calculate the series
(\ref{eq.(2.17)})-(\ref{eq.(2.20)}) we use the following
well-known formula from the residual theory (see [6], p.296)
\begin{equation}\label{eq.(2.21)}
\sum\limits_{\beta=-\infty}^{\infty}f(\beta)=-\sum\limits_{z_1,z_2,...,z_n}
\mathrm{res}(\pi \cot(\pi z)\cdot f(z)),
\end{equation}
where $z_1,z_2,...,z_n$ are poles of the function $f(z)$.

At first we consider $S_1$. We denote
$f_1(z)=\frac{1}{(z-h(p+\frac{\omega}{2\pi}))(z-h(p-\frac{\omega}{2\pi}))}$.
It is clear that $z_1=h(p+\frac{\omega}{2\pi})$ and
$z_2=h(p-\frac{\omega}{2\pi})$ are the poles of order 1 of the
function $f_1(z)$. Then taking into account the formula
(\ref{eq.(2.21)}) we have
\begin{equation}\label{eq.(2.22)}
S_1=\sum\limits_{\beta=-\infty}^{\infty}f_1(\beta)=-\sum\limits_{z_1,z_2}\mathrm{res}(\pi
\cot(\pi z)\cdot f_1(z)).
\end{equation}
Since
$$
\mathop{\mbox{res}}\limits_{z=z_1}(\pi \cot(\pi z)\cdot
f_1(z))=\lim\limits_{z\to z_1}\frac{\pi\cot(\pi
z)}{z-h(p-\frac{\omega}{2\pi})}=\frac{\pi^2}{h\omega}\cot(\pi
hp+\frac{h\omega}{2}),
$$
and
$$
\mathop{\mbox{res}}\limits_{z=z_2}(\pi \cot(\pi z)\cdot
f_1(z))=\lim\limits_{z\to z_2}\frac{\pi\cot(\pi
z)}{z-h(p+\frac{\omega}{2\pi})}=-\frac{\pi^2}{h\omega}\cot(\pi
hp-\frac{h\omega}{2}).
$$
Using the last two equalities, from (\ref{eq.(2.22)}) we get
$$
S_1=\frac{\pi^2}{h\omega}\left(\cot(\pi
hp-\frac{h\omega}{2})-\cot(\pi hp+\frac{h\omega}{2})\right)
=\frac{\pi^2\sin(h\omega)}{h\omega\sin(\pi
hp-\frac{h\omega}{2})\sin(\pi hp+\frac{h\omega}{2})}.
$$

Taking into account that $\lambda=e^{2\pi iph}$ and using the well
known formulas
$$
\cos z=\frac{e^{zi}+e^{-zi}}{2},\ \ \sin
z=\frac{e^{zi}-e^{-zi}}{2i},
$$
after some simplifications for $S_1$ we have
$$
S_1=\frac{-(2\pi)^2\lambda\sin h\omega}{h\omega(\lambda^2-2\lambda
\cos h\omega+1)}.
$$
Similarly for $S_2$ and $S_3$ we arrive the following equalities:
$$
S_2=\frac{-(2\pi)^2\lambda}{(\lambda^2+1)\cos
h\omega-2\lambda+i(\lambda^2-1)\sin h\omega}
$$
and
$$
S_3=\frac{-(2\pi)^2\lambda}{(\lambda^2+1)\cos
h\omega-2\lambda-i(\lambda^2-1)\sin h\omega}.
$$

The series (\ref{eq.(2.20)}) was calculated in \cite{Shad85},
where the following expression has been obtained
$$
F[h\stackrel{\abd}{G}_{k,1}]=\frac{h^{2k}\lambda
E_{2k-2}(\lambda)}{(2k-1)!\ (1-\lambda)^{2k}}.
$$
So, Lemma \ref{L2.4} is proved. \hfill $\Box$

\section{Construction of the discrete analogue of the operator}
\setcounter{equation}{0}

In this section we construct the discrete analogue $D_m(h\beta)$
of the operator $\frac{\d^{2m}}{\d
x^{2m}}+2\omega^2\frac{\d^{2m-2}}{\d
x^{2m-2}}+\omega^4\frac{\d^{2m-4}}{\d x^{2m-4}}$ and obtain some
properties of the constructed discrete function $D_m(h\beta)$.

The result is the following one:

\begin{theorem}\label{THM3.1}
 The discrete analogue to the differential operator
$\frac{\d^{2m}}{\d x^{2m}}+2\omega^2\frac{\d^{2m-2}}{\d
x^{2m-2}}+\omega^4\frac{\d^{2m-4}}{\d x^{2m-4}}$ satisfying
equation \emph{(\ref{eq.(1.1)})} has the form
\begin{equation}\label{eq.(3.1)}
D_m(h\beta)=\frac{2\omega^{2m-1}}{(-1)^mp_{2m-2}^{(2m-2)}}\left\{
\begin{array}{ll}
\sum\limits_{k=1}^{m-1}A_k \lambda_k^{|\beta|-1},& |\beta|\geq 2,\\
1+\sum\limits_{k=1}^{m-1}A_k,& |\beta|=1,\\
C+\sum\limits_{k=1}^{m-1}\frac{A_k}{\lambda_k},& \beta=0,
\end{array}
\right.
\end{equation}
where
\begin{eqnarray}
A_k&=&\frac{(1-\lambda_k)^{2m-4}(\lambda_k^2-2\lambda_k\cos
h\omega+1)^2p_{2m-2}^{(2m-2)}}{\lambda_k{\cal P}_{2m-2}'(\lambda_k)},\label{eq.(3.2)}\\
C&=&4-4\cos h\omega-2m-\frac{p_{2m-3}^{(2m-2)}}{p_{2m-2}^{(2m-2)}},\label{eq.(3.3)}\\
p_{2m-2}^{(2m-2)}&=&(2m-3)\sin h\omega-h\omega\cos h\omega\nonumber \\
&&\qquad\qquad\qquad
+2\sum\limits_{k=1}^{m-2}\frac{(-1)^k(m-k-1)(h\omega)^{2k-1}}{(2k-1)!},\label{eq.(3.4)}
\end{eqnarray}
${\cal P}_{2m-2}(\lambda)$ is the polynomial of degree $2m-2$
defined by \emph{(\ref{eq.(2.13)})}, $p_{2m-2}^{(2m-2)},\
p_{2m-3}^{(2m-2)}$ are the coefficients and $\lambda_k$ are the
roots of the polynomial ${\cal P}_{2m-2}(\lambda)$,
$|\lambda_k|<1$.
\end{theorem}

Before proving this result we present the following properties
that are not difficult to verify.

\begin{lemma}
The discrete analogue  $D_m(h\beta)$ of the differential operator
$\frac{\d^{2m}}{\d x^{2m}}+2\omega^2\frac{\d^{2m-2}}{\d
x^{2m-2}}+\omega^4\frac{\d^{2m-4}}{\d x^{2m-4}}$ satisfies the
following equalities

\emph{1)} $D_m(h\beta)*\sin(h\omega\beta)=0,$

\emph{2)} $D_m(h\beta)*\cos(h\omega\beta)=0,$

\emph{3)} $D_m(h\beta)*(h\omega\beta)\sin(h\omega\beta)=0,$

\emph{4)} $D_m(h\beta)*(h\omega\beta)\cos(h\omega\beta)=0,$

\emph{5)} $D_m(h\beta)*(h\beta)^\alpha=0,$ $\alpha=0,1,...,2m-5$.
\end{lemma}

\emph{Proof of Theorem \ref{THM3.1}.} According to the theory of
periodic distribution (generalized) functions and Fourier
transformations, instead of the function $D_m(h\beta)$ it is
convenient to search the harrow-shaped function (see
\cite{Sob74,SobVas})
$$
\stackrel{\abd}{D}_m\!\!(x)=\sum\limits_{\beta=-\infty}^{\infty}
D_m(h\beta)\delta(x-h\beta).
$$
In the class of harrow-shaped functions, equation (\ref{eq.(1.1)})
takes the following form
\begin{equation}\label{eq.(3.6)}
\stackrel{\abd}{D}_m\!\!(x)*\stackrel{\abd}{G}_m\!\!(x)=\delta(x),
\end{equation}
where
$\stackrel{\abd}{G}_m\!\!(x)=\sum\limits_{\beta=-\infty}^{\infty}
G_m(h\beta)\delta(x-h\beta)$ is the harrow-shaped function
corresponding to the discrete function $G_m(h\beta)$.

Applying the Fourier transformation to both sides of equation
(\ref{eq.(3.6)}) and taking into account (\ref{eq.(2.2)}) and
(\ref{eq.(2.4)}) we have
\begin{equation}\label{eq.(3.7)}
F[\stackrel{\abd}{D}_m\!\!(x)]=\frac{1}{F[\stackrel{\abd}{G}_m\!\!(x)]}.
\end{equation}
First, we calculate the Fourier transformation
$F[\stackrel{\abd}{G}_m\!\!(x)]$ of the harrow-shaped function
$\stackrel{\abd}{G}_m\!\!(x)$.

Using equalities (\ref{eq.(2.5)}), (\ref{eq.(2.6)}) and
(\ref{eq.(2.8)}), we get
\begin{eqnarray*}
\stackrel{\abd}{G}_m\!\!(x)&=&\sum_{\beta=-\infty}^{\infty}G_m(h\beta)\delta(x-h\beta)\nonumber\\
&=&h^{-1}G_m(x)\sum_{\beta=-\infty}^{\infty}\delta(h^{-1}x-\beta).\nonumber
\end{eqnarray*}
Hence
\begin{equation}
\stackrel{\abd}{G}_m\!\!(x)=h^{-1}G_m(x)\cdot
\phi_0(h^{-1}x).\label{eq.(3.8)}
\end{equation}
Taking into account (\ref{eq.(2.1)}) and (\ref{eq.(2.8)}), for
$F[\phi_0(h^{-1}x)]$ we obtain
\begin{eqnarray}
F[\phi_0(h^{-1}x)]&=&\sum\limits_{\beta=-\infty}^{\infty}\int_{-\infty}^{\infty}\delta(h^{-1}x-\beta)e^{2\pi ipx}\d x\nonumber\\
&=&h\sum\limits_{\beta=-\infty}^{\infty}\int_{-\infty}^{\infty}\delta(x-h\beta)e^{2\pi ipx}\d x\nonumber\\
&=&h\sum\limits_{\beta=-\infty}^{\infty}e^{2\pi
iph\beta}=h\sum\limits_{\beta=-\infty}^{\infty}\delta(hp-\beta)=h\phi_0(hp).\label{eq.(3.9)}
\end{eqnarray}
Applying the Fourier transformation to both sides of
(\ref{eq.(3.8)}) and using (\ref{eq.(2.3)}), (\ref{eq.(3.9)}), we
get
\begin{equation}
F[\stackrel{\abd}{G}_m\!\!(x)]=F[G_m(x)]*\phi_0(hp).
\label{eq.(3.10)}
\end{equation}
To calculate the Fourier transformation $F[G_m(x)]$, we use
equation (\ref{eq.(1.4)}). Taking into account (\ref{eq.(2.7)}),
we rewrite equation (\ref{eq.(1.4)}) in the following form
$$
(\delta^{(4)}(x)+2\omega^2\delta^{(2m-2)}(x)+\omega^4\delta^{(2m-4)}(x))*G_m(x)=\delta(x).
$$
Hence, keeping in mind  (\ref{eq.(2.2)}) and (\ref{eq.(2.4)}), we
have
$$
F[G_m(x)]=\frac{1}{(2\pi ip)^{2m}+2\omega^2(2\pi
ip)^{2m-2}+\omega^4 (2\pi ip)^{2m-4}}.
$$
Using the last equality, from (\ref{eq.(3.10)}) we arrive to
\begin{eqnarray*}
F[\stackrel{\abd}{G}_m\!\!(x)]&=&\frac{1}{(2\pi
ip)^{2m}+2\omega^2(2\pi ip)^{2m-2}+\omega^4 (2\pi
ip)^{2m-4}}*\phi_0(hp)\\
&=&
\frac{1}{h}\sum\limits_{\beta=-\infty}^{\infty}\int\limits_{-\infty}^{\infty}\frac{\delta(y-h^{-1}\beta)\d
y}{(2\pi i(p-y))^{2m}+2\omega^2(2\pi i(p-y))^{2m-2}+\omega^4 (2\pi
i(p-y))^{2m-4}}.
\end{eqnarray*}
Hence, after some calculations we get
\begin{equation}\label{eq.(3.11)}
F[\stackrel{\abd}{G}_m\!\!(x)]=
\sum\limits_{\beta=-\infty}^{\infty}\frac{h^{-1}}{(2\pi
i(p-\frac{\beta}{h}))^{2m}+2\omega^2(2\pi
i(p-\frac{\beta}{h}))^{2m-2}+\omega^4 (2\pi
i(p-\frac{\beta}{h}))^{2m-4}}.
\end{equation}
Now, expanding to partial fractions of the right hand side of
(\ref{eq.(3.11)}) and taking into account (\ref{eq.(3.7)}), we
conclude
\begin{eqnarray}
F[\stackrel{\abd}{D}_m](p)&=&h\Bigg[\frac{(-1)^{m-1}(2m-3)h^2}{(2\pi)^2\cdot
2\omega^{2m-2}}S_1+\frac{(-1)^{m-2}h^2}{(2\pi)^2\cdot
4\omega^{2m-2}}(S_2+S_3) \nonumber\\
&&\qquad+\frac{1}{\omega^{2m}}\sum_{k=1}^{m-2}(-1)^{m-k}(m-k-1)\omega^{2k}F[h\stackrel{\abd}{G}_{k,1}]\Bigg]^{-1},\label{eq.(3.12)}
\end{eqnarray}
where $S_1$, $S_2$, $S_3$ and $F[h\stackrel{\abd}{G}_{k,1}]$ are
defined by (\ref{eq.(2.17)})-(\ref{eq.(2.20)}).

It is clear that the function $F[\stackrel{\abd}{D}_m](p)$ is a
periodic function with respect to the variable $p$, and the period
is $h^{-1}$. Moreover it is real and analytic for all real $p$.
The zeros of the function $F[\stackrel{\abd}{D}_m](p)$ are
$p=h^{-1}\beta+\frac{\omega}{2\pi}$,
$p=h^{-1}\beta-\frac{\omega}{2\pi}$ and $p=h^{-1}\beta$.

The function $F[\stackrel{\abd}{D}_m ](p)$ can be expanded to the
Fourier series as follows
\begin{equation}
F[\stackrel{\abd}{D}_m](p)=\sum\limits_{\beta=-\infty}^{\infty}\hat
D_m(h\beta)e^{2\pi i ph\beta},\label{eq.(3.12)}
\end{equation}
where $\hat D_m(h\beta)$ are the Fourier coefficients of the
function $F[\stackrel{\abd}{D}_m](p)$, i.e.
\begin{equation}\label{eq.(3.13)}
\hat D_m(h\beta)=\int\limits_0^{h^{-1}}F[\stackrel{\abd}{D}_m](p)\
e^{-2\pi i ph\beta}\d p.
\end{equation}
Applying the inverse Fourier transformation to both sides of
(\ref{eq.(3.12)}), we obtain the following harrow-shaped function
$$
\stackrel{\abd}{D}_m\!\!(x)=\sum\limits_{\beta=-\infty}^{\infty}\hat
D_m(h\beta)\delta(x-h\beta).
$$
Then, according to the definition of harrow-shaped functions, the
discrete function $\hat D_m(h\beta)$ is the discrete argument
function $D_m(h\beta)$, which we are searching. So, we have to
find the discrete argument function $\hat D_m(h\beta)$.

But here, to find the function $\hat D_m(h\beta)$ we will not use
the formula (\ref{eq.(3.13)}). In this case to obtain $\hat
D_m(h\beta)$ we have to calculate the series $S_1$, $S_2$, $S_3$
and $F[h\stackrel{\abd}{G}_{k,1}]$ defined by
(\ref{eq.(2.17)})-(\ref{eq.(2.20)}). We have calculated these
series in Lemma \ref{L2.4}. Therefore, using Lemma \ref{L2.4},
from (\ref{eq.(3.13)}) we get
\begin{equation}
F[\stackrel{\abd}{D}_m](p)=\frac{2\omega^{2m-1}}{(-1)^m}\
\frac{(1-\lambda)^{2m-4}(\lambda^2+1-2\lambda\cos
h\omega)^2}{\lambda\ {\cal P}_{2m-2}(\lambda)},
\end{equation}
where $\lambda=e^{2\pi iph}$ and ${\cal P}_{2m-2}(\lambda)$ is the
polynomial of degree $2m-2$ defined by (\ref{eq.(2.13)}).

Now we will get the Fourier series of the function
$F[\stackrel{\abd}{D}_m](p)$. Dividing the polynomial
${(1-\lambda)^{2m-4}(\lambda^2+1-2\lambda\cos h\omega)^2}$ by the
polynomial $\lambda {\cal P}_{2m-2}(\lambda)$ we obtain
\begin{eqnarray}\label{eq.(3.15)}
&&F[\stackrel{\abd}{D}_m](p)=\frac{2\omega^{2m-1}}{(-1)^m}\
\frac{(1-\lambda)^{2m-4}(\lambda^2+1-2\lambda\cos
h\omega)^2}{\lambda\ {\cal P}_{2m-2}(\lambda)}\nonumber\\
&&
=\frac{2\omega^{2m-1}}{(-1)^mp_{2m-2}^{(2m-2)}}\left[\lambda+4-4\cos
h\omega-2m-\frac{p_{2m-3}^{(2m-2)}}{p_{2m-2}^{(2m-2)}}\right]+\frac{S_{2m-2}(\lambda)}{\lambda{\cal
P}_{2m-2}(\lambda)},
\end{eqnarray}
where $S_{2m-2}(\lambda)$ is a polynomial of degree $2m-2$.

It is clear that $\frac{S_{2m-2}(\lambda)}{\lambda{\cal
P}_{2m-2}(\lambda)}$  is the proper fraction. Since the roots of
the polynomial ${\lambda{\cal P}_{2m-2}(\lambda)}$ are real and
simple (see Appendix), then the rational fraction
$\frac{S_{2m-2}(\lambda)}{\lambda{\cal P}_{2m-2}(\lambda)}$ is
expanded to the sum of partial fractions, i.e.
\begin{equation}\label{eq.(3.16)}
\frac{S_{2m-2}(\lambda)}{\lambda {\cal
P}_{2m-2}(\lambda)}=\frac{1}{p_{2m-2}^{(2m-2)}}\left[\frac{A_0}{\lambda}+\sum\limits_{k=1}^{m-1}\left(\frac{A_{1,k}}{\lambda-\lambda_{1,k}}+
\frac{A_{2,k}}{\lambda-\lambda_{2,k}}\right)\right],
\end{equation}
where $A_0,$ $A_{1,k}$ and $A_{2,k}$ are unknown coefficients,
$\lambda_{1,k}$, $\lambda_{2,k}$ are the roots of the polynomial
${\cal P}_{2m-2}(\lambda)$ and $|\lambda_{1,k}|<1$,
$|\lambda_{2,k}|>1$ and
\begin{equation}\label{eq.(3.17)}
\lambda_{1,k}\lambda_{2,k}=1.
\end{equation}

Keeping in mind equality (\ref{eq.(3.16)}) from (\ref{eq.(3.15)}),
we obtain
\begin{eqnarray}
&&F[\stackrel{\abd}{D}_m](p)=\frac{2\omega^{2m-1}}{(-1)^m}\
\frac{(1-\lambda)^{2m-4}(\lambda^2+1-2\lambda\cos
h\omega)^2}{\lambda\ {\cal P}_{2m-2}(\lambda)}\label{eq.(3.18)}\\
&& \qquad\qquad\qquad
=\frac{2\omega^{2m-1}}{(-1)^mp_{2m-2}^{(2m-2)}}\left[\lambda+4-4\cos
h\omega-2m-\frac{p_{2m-3}^{(2m-2)}}{p_{2m-2}^{(2m-2)}}\right]\nonumber\\
&&\qquad\qquad\qquad\ \ \
+\frac{1}{p_{2m-2}^{(2m-2)}}\left[\frac{A_0}{\lambda}+\sum\limits_{k=1}^{m-1}\left(\frac{A_{1,k}}{\lambda-\lambda_{1,k}}+
\frac{A_{2,k}}{\lambda-\lambda_{2,k}}\right)\right].\nonumber
\end{eqnarray}
Multiplying both sides of the last equality by $\lambda{\cal
P}_{2m-2}(\lambda)$ we have
\begin{eqnarray}
&&\frac{2\omega^{2m-1}}{(-1)^m}\
{(1-\lambda)^{2m-4}(\lambda^2+1-2\lambda\cos
h\omega)^2}\nonumber\\
&& \qquad =\lambda\ {\cal
P}_{2m-2}(\lambda)\Bigg[\frac{2\omega^{2m-1}}{(-1)^mp_{2m-2}^{(2m-2)}}\left[\lambda+4-4\cos
h\omega-2m-\frac{p_{2m-3}^{(2m-2)}}{p_{2m-2}^{(2m-2)}}\right]\nonumber\\
&&\qquad\qquad\qquad\ \ \
+\frac{1}{p_{2m-2}^{(2m-2)}}\left[\frac{A_0}{\lambda}+\sum\limits_{k=1}^{m-1}\left(\frac{A_{1,k}}{\lambda-\lambda_{1,k}}+
\frac{A_{2,k}}{\lambda-\lambda_{2,k}}\right)\right]\Bigg].\label{eq.(3.19)}
\end{eqnarray}
To find unknown coefficients $A_0$, $A_{1,k}$ and $A_{2,k}$ in
(\ref{eq.(3.19)}) we put $\lambda=0$, $\lambda=\lambda_{1,k}$ and
$\lambda=\lambda_{2,k}$. Then  we get
\begin{eqnarray}
A_0&=&\frac{2\omega^{2m-1}}{(-1)^m},\nonumber\\
A_{1,k}&=&\frac{2\omega^{2m-1}(1-\lambda_{1,k})^{2m-4}(\lambda_{1,k}^2+1-2\lambda_{1,k}\cos
h\omega)^2p_{2m-2}^{(2m-2)}}{(-1)^m\lambda_{1,k}{\cal
P}_{2m-2}'(\lambda_{1,k})},\label{eq.(3.20)}\\
A_{2,k}&=&\frac{2\omega^{2m-1}(1-\lambda_{2,k})^{2m-4}(\lambda_{2,k}^2+1-2\lambda_{2,k}\cos
h\omega)^2p_{2m-2}^{(2m-2)}}{(-1)^m\lambda_{2,k}{\cal
P}_{2m-2}'(\lambda_{2,k})}.\nonumber
\end{eqnarray}
Hence, taking into account (\ref{eq.(3.17)}) and Lemma \ref{L2.3},
we get
\begin{equation}\label{eq.(3.21)}
A_{2,k}=-\frac{A_{1,k}}{\lambda_{1,k}^2}.
\end{equation}
Since $|\lambda_{1,k}|<1$ and $|\lambda_{2,k}|>1$ the expressions
$$
\sum\limits_{k=1}^{m-1}\frac{A_{1,k}}{\lambda-\lambda_{1,k}}\mbox{
and } \sum\limits_{k=1}^{m-1}\frac{A_{2,k}}{\lambda-\lambda_{2,k}}
$$
can be expanded to the Laurent series on the circle $|\lambda|
=1$, i.e.
$$
\sum\limits_{k=1}^{m-1}\frac{A_{1,k}}{\lambda-\lambda_{1,k}}=\frac{1}{\lambda}
\sum\limits_{k=1}^{m-1}\frac{A_{1,k}}{1-\frac{\lambda_{1,k}}{\lambda}}=\frac{1}{\lambda}
\sum\limits_{k=1}^{m-1}A_{1,k}\sum\limits_{\beta=0}^{\infty}\left(\frac{\lambda_{1,k}}{\lambda}\right)^{\beta}
$$
and taking into account (\ref{eq.(3.17)}), we deduce
$$
\sum\limits_{k=1}^{m-1}\frac{A_{2,k}}{\lambda-\lambda_{2,k}}=
\sum\limits_{k=1}^{m-1}\frac{A_{2,k}}{\lambda-\frac{1}{\lambda_{1,k}}}=
-\sum\limits_{k=1}^{m-1}\frac{A_{2,k}\lambda_{1,k}}{1-\lambda\
\lambda_{1,k}}=-\sum_{k=1}^{m-1}A_{2,k}\lambda_{1,k}\sum\limits_{\beta=0}^{\infty}(\lambda\
\lambda_{1,k})^\beta.
$$
Using the last two equalities and taking into account
(\ref{eq.(3.3)}) and $\lambda=e^{2\pi iph}$ from (\ref{eq.(3.18)})
we get
\begin{eqnarray*}
F[\stackrel{\abd}{D}_m](p)&=&\frac{1}{p_{2m-2}^{(2m-2)}}\Bigg\{\frac{2\omega^{2m-1}}{(-1)^m}[e^{2\pi
iph}+C]+A_0e^{-2\pi iph}\\
&&+\sum\limits_{k=1}^{m-1}A_{1,k}\sum_{\beta=0}^{\infty}\lambda_{1,k}^\beta
e^{-2\pi
iph(\gamma+1)}-\sum_{k=1}^{m-1}A_{2,k}\sum_{\beta=0}^{\infty}\lambda_{1,k}^{\beta+1}e^{2\pi
iph\gamma}\Bigg\}.
\end{eqnarray*}
Thus the Fourier series for $F[\stackrel{\abd}{D}_m](p)$ has the
following form
$$
F[\stackrel{\abd}{D}_m](p)=\sum\limits_{\beta=-\infty}^{\infty}D_m(h\beta)e^{2\pi
i hp\beta},
$$
where
$$
D_m(h\beta)=\frac{1}{p_{2m-2}^{(2m-2)}}\left\{
\begin{array}{ll}
\sum\limits_{k=1}^{m-1}A_{1,k}\lambda_{1,k}^{-\beta-1},&\beta\leq
-2,\\
A_0+\sum\limits_{k=1}^{m-1}A_{1,k},& \beta=-1,\\
\frac{2\omega^{2m-1}}{(-1)^m}C-\sum\limits_{k=1}^{m-1}A_{2,k}\lambda_{1,k},&
\beta=0,\\
\frac{2\omega^{2m-1}}{(-1)^m}-\sum\limits_{k=1}^{m-1}A_{2,k}\lambda_{1,k}^2,&
\beta=1,\\
-\sum\limits_{k=1}^{m-1}A_{2,k}\lambda_{1,k}^{\beta+1},&\beta\geq
2.
\end{array}
\right.
$$
Thus, using (\ref{eq.(3.21)}), we rewrite the discrete function
$D_m(h\beta)$ in the following form
\begin{equation}\label{eq.(3.22)}
D_m(h\beta)=\frac{1}{p_{2m-2}^{(2m-2)}}\left\{
\begin{array}{ll}
\sum\limits_{k=1}^{m-1}A_{1,k}\lambda_{1,k}^{|\beta|-1},&|\beta|\geq
2,\\
\frac{2\omega^{2m-1}}{(-1)^m}+\sum\limits_{k=1}^{m-1}A_{1,k},& |\beta|=1,\\
\frac{2\omega^{2m-1}}{(-1)^m}C+\sum\limits_{k=1}^{m-1}\frac{A_{1,k}}{\lambda_{1,k}},&
\beta=0.\\
\end{array}
\right.
\end{equation}
Combining (\ref{eq.(3.2)}) and (\ref{eq.(3.20)}) we get
\begin{equation}\label{eq.(3.23)}
A_{1,k}=\frac{2\omega^{2m-1}}{(-1)^m}\ A_k.
\end{equation}
Finally, denoting $\lambda_k=\lambda_{1,k}$ and taking into
account (\ref{eq.(3.23)}) from (\ref{eq.(3.22)}) we get the
statement of the Theorem. Theorem \ref{THM3.1} is proved.\hfill
$\Box$

\textit{Remark 1.} {From (\ref{eq.(3.1)}) we note that
$D_m(h\beta)$ is an even function, i.e.
$$D_m(h\beta)=D_m(-h\beta)$$ and since $|\lambda_k|<1$, then the
function $D_m(h\beta)$ decreases exponentially as $|\beta|\to
\infty$.}

\textit{Remark 2.} It should be noted that as consequence Theorem
\ref{THM3.1} generalizes previous results given in
\cite{Hay09,Hay12} for the particular cases $m=2,3$.

\section{Appendix}
\setcounter{equation}{0}

Here we obtain explicate formulas for the coefficients of the
polynomials ${\cal P}_{2m-2}(x)$ defined by (\ref{eq.(2.13)}), for
the cases $m=2,3,...,6$. Furthermore, for these cases, we show
connection of the polynomials ${\cal P}_{2m-2}(x)$ with the
Euler-Frobenius polynomials $E_{2m-2}(x)$ defined by
(\ref{eq.(2.9)}).

Using equality (\ref{eq.(2.11)}), when $m=2,3,...,6$, for the
Euler-Frobenius polynomials $E_{2m-2}(x)$ we have
\begin{eqnarray}\label{eq.(4.1)}
E_2(x)&=&x^2+4x+1,\nonumber\\
E_4(x)&=&x^4+26x^3+66x^2+26x+1,\nonumber\\
E_6(x)&=&x^6+120x^5+1191x^4+2416x^3+1191x^2+120x+1, \nonumber\\
E_8(x)&=&x^8+502x^7+14608x^6+88234x^5+156190x^4+88234x^3\\
&&+14608x^2+502x+1, \nonumber\\
E_{10}(x)&=&x^{10}+2036x^9+152637x^8+2203488x^7+9738114x^6\nonumber\\
&&+15724248x^5+9738114x^4+2203488x^3+152637x^2+2036x+1. \nonumber
\end{eqnarray}

Now, by direct computation from (\ref{eq.(2.13)}) when
$m=2,3,...,6$, for the coefficients $p_s^{(2m-2)}$,
$s=0,1,...,2m-2$, of the polynomials ${\cal P}_{2m-2}(x)$ we
consequently get the following results.

1) For $m=2$: ${\cal P}_2(x)=\sum\limits_{s=0}^2p_s^{(2)}$ and
\begin{eqnarray}\label{eq.(4.2)}
p_2^{(2)}&=&p_0^{(2)}=\sin(h\omega)-h\omega\cos(h\omega),\nonumber\\
p_1^{(2)}&=&2h\omega-\sin(2h\omega),\nonumber\\
\frac{p_2^{(2)}}{p_2^{(2)}}&=&\frac{p_0^{(2)}}{p_2^{(2)}}=1,\\
\frac{p_1^{(2)}}{p_2^{(2)}}&=&4-\frac25(h\omega)^2+O(h^4).\nonumber
\end{eqnarray}

2) For $m=3$: ${\cal P}_4(x)=\sum\limits_{s=0}^4p_s^{(4)}$ and
\begin{eqnarray}\label{eq.(4.3)}
p_4^{(4)}&=&p_0^{(4)}=3\sin(h\omega)-h\omega\cos(h\omega)-2h\omega,\nonumber\\
p_3^{(4)}&=&p_1^{(4)}=10h\omega\cos(h\omega)+2h\omega-3\sin(2h\omega)-6\sin(h\omega),\nonumber\\
p_2^{(4)}&=&-8h\omega(1+\cos^2(h\omega))+6\sin(h\omega)-2h\omega\cos(h\omega)+6\sin(2h\omega),\nonumber\\
\frac{p_4^{(4)}}{p_4^{(4)}}&=&\frac{p_0^{(4)}}{p_4^{(4)}}=1,\\
\frac{p_3^{(4)}}{p_4^{(4)}}&=&\frac{p_1^{(4)}}{p_4^{(4)}}=26-\frac{18}{7}(h\omega)^2+O(h^4),\nonumber\\
\frac{p_2^{(4)}}{p_4^{(4)}}&=&66-\frac{64}{7}(h\omega)^2+O(h^4).\nonumber
\end{eqnarray}

3) For $m=4$: ${\cal P}_6(x)=\sum\limits_{s=0}^6p_s^{(6)}$ and
\begin{eqnarray}\label{eq.(4.4)}
p_6^{(6)}&=&p_0^{(6)}=5\sin(h\omega)-h\omega\cos(h\omega)-4h\omega+\frac13(h\omega)^3,\nonumber\\
p_5^{(6)}&=&p_1^{(6)}=10h\omega-5\sin(2h\omega)-20\sin(h\omega)+4h\omega\cos(h\omega)+\frac43(h\omega)^3\nonumber\\
&&\ \ \ \ \ \ \ -8\cos(h\omega)\left(-2h\omega+\frac16(h\omega)^3\right),\nonumber\\
p_4^{(6)}&=&p_2^{(6)}=35\sin(h\omega)-7h\omega\cos(h\omega)-12h\omega+20\sin(2h\omega)+\frac13(h\omega)^3\nonumber\\
&&\ \ \ \ \ \ \ -8\cos(h\omega)\left(4h\omega+\frac23(h\omega)^3\right)+4(1+2\cos^2(h\omega))\left(-2h\omega+\frac16(h\omega)^3\right),\nonumber\\
p_3^{(6)}&=&-16\cos(h\omega)\left(-2h\omega+\frac16(h\omega)^3\right)+4(1+2\cos^2(h\omega))\left(4h\omega+\frac23(h\omega)^3\right)\nonumber\\
&&-40\sin(h\omega)+8h\omega\cos(h\omega)+12h\omega-30\sin(2h\omega),\nonumber\\
\frac{p_6^{(6)}}{p_6^{(6)}}&=&\frac{p_0^{(6)}}{p_6^{(6)}}=1,\\
\frac{p_5^{(6)}}{p_6^{(6)}}&=&\frac{p_1^{(6)}}{p_6^{(6)}}=120-\frac{26}{3}(h\omega)^2+O(h^4),\nonumber\\
\frac{p_4^{(6)}}{p_6^{(6)}}&=&\frac{p_2^{(6)}}{p_6^{(6)}}=1191-\frac{470}{3}(h\omega)^2+O(h^4),\nonumber\\
\frac{p_3^{(6)}}{p_6^{(6)}}&=&2416-\frac{1108}{3}(h\omega)^2+O(h^4).\nonumber
\end{eqnarray}

4) For $m=5$: ${\cal P}_8(x)=\sum\limits_{s=0}^8p_s^{(8)}$ and
\begin{eqnarray}\label{eq.(4.5)}
p_8^{(8)}&=&p_0^{(8)}=7\sin(h\omega)-h\omega\cos(h\omega)-6h\omega+\frac{2}{3}(h\omega)^3-\frac{1}{60}(h\omega)^5,\nonumber\\
p_7^{(8)}&=&p_1^{(8)}=-\frac{13}{30}(h\omega)^5+26h\omega+\frac{4}{3}(hw)^3-8\cos(h\omega)\left(-3h\omega+\frac{1}{3}(h\omega)^3-\frac{1}{120}(h\omega)^5\right)\nonumber\\
&&\ \ \ \ \ \ \ -7\sin(2h\omega)-42\sin(h\omega)+6h\omega\cos(h\omega),\nonumber\\
p_6^{(8)}&=&p_2^{(8)}=-\frac{11}{10}(h\omega)^5-4(h\omega)^3-48h\omega-8\cos(h\omega)\left(-\frac{13}{60}(h\omega)^5+12h\omega+\frac{2}{3}(h\omega)^3\right)\nonumber\\
& &\ \ \ \ \ \ \  +4(1+2\cos^2(h\omega))\left(-3h\omega+\frac{1}{3}(h\omega)^3-\frac{1}{120}(h\omega)^5\right)+112\sin(h\omega)\nonumber\\
& &\ \ \ \ \ \ \ -16h\omega\cos(h\omega)+42\sin(2h\omega),\nonumber\\
p_5^{(8)}&=&p_3^{(8)}=-\frac{13}{30}(h\omega)^5+54h\omega+\frac{4}{3}(h\omega)^3+8\cos(h\omega)\left(\frac{67}{120}(h\omega)^5+\frac{5}{3}(h\omega)^3+21h\omega\right)\nonumber\\
&&\ \ \ \ \ \ \ +4(1+2\cos^2(h\omega))\left(-\frac{13}{60}(h\omega)^5+12h\omega+\frac{2}{3}(h\omega)^3\right)-182\sin(h\omega)\nonumber\\
&&\ \ \ \ \ \ \ +26h\omega\cos(h\omega)-105\sin(2h\omega),\nonumber\\
p_4^{(8)}&=&210\sin(h\omega)-30h\omega\cos(h\omega)-52h\omega+140\sin(2h\omega)+\frac{4}{3}(h\omega)^3-\frac{1}{30}(h\omega)^5\nonumber\\
&&-16\cos(h\omega)
\left(-\frac{13}{60}(h\omega)^5+12h\omega+\frac{2}{3}(h\omega)^3\right)\nonumber\\
&&+4(1+2\cos^2(h\omega))\left(-\frac{11}{20}(h\omega)^5-2(h\omega)^3-18h\omega\right),\nonumber\\
\frac{p_8^{(8)}}{p_8^{(8)}}&=&\frac{p_0^{(8)}}{p_8^{(8)}}=1,\\
\frac{p_7^{(8)}}{p_8^{(8)}}&=&\frac{p_1^{(8)}}{p_8^{(8)}}=502-\frac{1426}{55}(h\omega)^2+O(h^4),\nonumber\\
\frac{p_6^{(8)}}{p_8^{(8)}}&=&\frac{p_2^{(8)}}{p_8^{(8)}}=14608-\frac{86664}{55}(h\omega)^2+O(h^4),\nonumber\\
\frac{p_5^{(8)}}{p_8^{(8)}}&=&\frac{p_3^{(8)}}{p_8^{(8)}}=88234-\frac{141798}{11}(h\omega)^2+O(h^4),\nonumber\\
\frac{p_4^{(8)}}{p_8^{(8)}}&=&156190-\frac{273872}{11}(h\omega)^2+O(h^4).\nonumber
\end{eqnarray}

5) For $m=6$: ${\cal P}_{10}(x)=\sum\limits_{s=0}^{10}p_s^{(10)}$
and
\begin{eqnarray}\label{eq.(4.6)}
p_{10}^{(10)}&=&p_0^{(10)}= 9\sin(h\omega)-h\omega\cos(h\omega)-8h\omega+\frac{1}{2520}(h\omega)^7+(hw)^3-\frac{1}{30}(h\omega)^5,\nonumber\\
p_9^{(10)}&=&p_1^{(10)}=50h\omega-9\sin(2h\omega)-72\sin(h\omega)+8h\omega\cos(h\omega)-\frac{4}{5}(h\omega)^5+\frac{1}{21}(h\omega)^7\nonumber\\
&&\ \ \ \ \ \ \ \  -8\cos(h\omega)\left(-4h\omega+\frac{1}{5040}(h\omega)^7+\frac{1}{2}(h\omega)^3-\frac{1}{60}(h\omega)^5\right),\nonumber\\
p_8^{(10)}&=&p_2^{(10)}=261\sin(h\omega)-29h\omega\cos(h\omega)-136h\omega+72\sin(2h\omega)+\frac{397}{840}(h\omega)^7-\frac{1}{2}(h\omega)^5\nonumber\\
&&\ \ \ \ \ \ \ \  -9(h\omega)^3-8\cos(h\omega)\left(24h\omega-\frac{2}{5}(h\omega)^5+\frac{1}{42}(h\omega)^7\right)\nonumber\\
&&\ \ \ \ \ \ \ \  +4(1+2\cos^2(h\omega))\left(-4h\omega+\frac{1}{5040}(h\omega)^7+\frac{1}{2}(h\omega)^3-\frac{1}{60}(h\omega)^5\right),    \nonumber\\
p_7^{(10)}&=&p_3^{(10)}=8\cos(h\omega)\left(64h\omega-\frac{149}{630}(h\omega)^7-\frac{4}{15}(h\omega)^5-4(h\omega)^3\right)+216h\omega+\frac{8}{3}(h\omega)^5\nonumber\\
&&\ \ \ \ \ \ \ \ +16(h\omega)^3+\frac{302}{315}(h\omega)^7+4(1+2\cos^2(h\omega))\left(24h\omega-\frac{2}{5}(h\omega)^5+\frac{1}{42}(h\omega)^7\right)\nonumber\\
&&\ \ \ \ \ \ \ \ -576\sin(h\omega)+64h\omega\cos(h\omega)-252\sin(2h\omega),\nonumber\\
p_6^{(10)}&=&p_4^{(10)}=882\sin(h\omega)-98h\omega\cos(h\omega)-240h\omega+504\sin(2h\omega)+\frac{149}{315}(h\omega)^7-\frac{8}{15}(h\omega)^5\nonumber\\
&&\ \ \ \ \ \ \ \
-8(h\omega)^3-8\cos(h\omega)\left(104h\omega+\frac{14}{15}(h\omega)^5+8(h\omega)^3+\frac{151}{315}(h\omega)^7\right)\nonumber\\
&&\ \ \ \ \ \ \ \ +4(1+2\cos^2(h\omega))\left(-60h\omega+\frac{397}{1680}(h\omega)^7-\frac{1}{4}(h\omega)^5-\frac{9}{2}(h\omega)^3\right),\nonumber\\
p_5^{(10)}&=&-1008\sin(h\omega)+112h\omega\cos(h\omega)+236h\omega-630\sin(2h\omega)-\frac{8}{5}(h\omega)^5+\frac{2}{21}(h\omega)^7\nonumber \\
&&\ \ \ \ \ \ \ \ -16\cos(h\omega)\left(-60h\omega+\frac{397}{1680}(h\omega)^7-\frac{1}{4}(h\omega)^5-\frac{9}{2}(h\omega)^3\right)\nonumber\\
&&\ \ \ \ \ \ \ \ +4(1+2\cos^2(h\omega))\left(80h\omega+\frac{4}{3}(h\omega)^5+8(h\omega)^3+\frac{151}{315}(h\omega)^7\right),\nonumber\\
\frac{p_{10}^{(10)}}{p_{10}^{(10)}}&=&\frac{p_0^{(10)}}{p_{10}^{(10)}}=1,\\
\frac{p_9^{(10)}}{p_{10}^{(10)}}&=&\frac{p_1^{(10)}}{p_{10}^{(10)}}=2036-\frac{998}{13}(h\omega)^2+O(h^4),\nonumber\\
\frac{p_8^{(10)}}{p_{10}^{(10)}}&=&\frac{p_2^{(10)}}{p_{10}^{(10)}}=152637-\frac{170920}{13}(h\omega)^2+O(h^4),\nonumber\\
\frac{p_7^{(10)}}{p_{10}^{(10)}}&=&\frac{p_3^{(10)}}{p_{10}^{(10)}}=2203488-\frac{3644480}{13}(h\omega)^2+O(h^4),\nonumber\\
\frac{p_6^{(10)}}{p_{10}^{(10)}}&=&\frac{p_4^{(10)}}{p_{10}^{(10)}}=9738114-\frac{19460312}{13}(h\omega)^2+O(h^4),\nonumber\\
\frac{p_5^{(10)}}{p_{10}^{(10)}}&=&15724248-\frac{33280180}{13}(h\omega)^2+O(h^4).\nonumber
\end{eqnarray}

From (\ref{eq.(4.1)})-(\ref{eq.(4.6)}) we state the following.

\begin{conjecture}\label{L2.2}
For the coefficients $p_s^{(2m-2)}$, $s=0,1,...,2m-2$ of the
polynomial ${\cal P}_{2m-2}(x)$ the following holds
$$
\frac{p_s^{(2m-2)}}{p_{2m-2}^{(2m-2)}}=a_s^{(2m-2)}+O(h^2), \ \ \
s=0,1,...,2m-2,
$$
where $a_s^{(2m-2)}$ are the coefficients of the Euler-Frobenius
polynomial $E_{2m-2}(x)$ of degree $2m-2$.
\end{conjecture}

\section*{Acknowledgements}
The author thanks professor A.Cabada for discussion of the
results. The present work was done in the University of Santiago
de Compostela, Spain. A.R. Hayotov thanks the program Erasmus
Mundus Action 2, Lot 10, Marco XXI for financial support (project
number: Lot 10 - 20112572).

\end{document}